\documentclass[11pt]{article}
\usepackage{mathrsfs}
\usepackage{amssymb}
\usepackage{amsfonts}
\usepackage{amsmath}
\usepackage{amscd}
\usepackage{amsthm}
\usepackage{latexsym}
\usepackage[all]{xy}
\usepackage{CJK}
\numberwithin{equation}{section}

\newcommand {\add}{\mathrm{add}}
\newcommand {\End}{\mathrm{End}}
\newcommand {\Mor}{\mathrm{Mor}}
\newcommand {\Hom}{\mathrm{Hom}}
\newcommand {\Ext}{\mathrm{Ext}_A^1}
\newcommand {\Mod}{\mathrm{mod}}
\newcommand {\ind}{\mathrm{ind}}
\newcommand {\module}{\textrm{module}}
\newcommand {\modules}{\textrm{modules}}

\newcommand {\mf}{\mathcal{F}}

\newcommand {\mc}{\mathcal{C}}
\newcommand {\gld}{\mathrm{gld}}

\newcommand {\Ker}{\mathrm{Ker}}

\newcommand{\ol}{\overline}

\newcommand {\dt}{DT}
\newcommand{\Tau}{\tau^{-1}}

\begin{document}
\title{\bf Representation dimensions of triangular matrix algebras$^\star$}
\author{{\small  Hongbo Yin, Shunhua Zhang$^*$}\\
{\small  School of Mathematics,\ Shandong University,\ Jinan 250100,
P. R. China }}

\pagenumbering{arabic}
\date{}
\maketitle
\begin{abstract}
\begin{center}
 \begin{minipage}{120mm}
   \small\rm
 {\bf  Abstract}\ \ Let $A$ be a finite dimensional hereditary algebra over an
algebraically closed field $k$, $T_2(A)=\left(\begin{array}{cc}A&0\\
A&A\end{array}\right)$ be the triangular matrix algebra and
$A^{(1)}=\left(\begin{array}{cc}A&0\\ DA&A\end{array}\right)$ be the
duplicated algebra of $A$ respectively. We prove that ${\rm
rep.dim}\ T_2(A)$ is at most three if $A$ is Dynkin type and ${\rm
rep.dim}\ T_2(A)$ is at most four if $A$ is not Dynkin type. Let $T$
be a  tilting A-$\module$ and $\ol{T}=T\oplus\ol{P}$ be a tilting
$A^{(1)}$-$\module$. We show that $\End_{A^{(1)}}~\ol{T}$ is
representation finite if and only if the full subcategory
$\{(X,Y,f)\ |\ X\in {\rm mod}\ A,
Y\in\tau^{-1}\mathscr{F}(T_A)\cup{\rm add}\ A\}$ of ${\rm mod \
T_2(A)}$ is of finite type, where $\tau$ is the
Auslander-Reiten translation and $\mathscr{F}(T_A)$ is the
torsion-free class of ${\rm mod}\ A$ associated with $T$. Moreover,
we also prove that ${\rm rep.dim\ End}_{A^{(1)}}\ {\ol T}$ is at
most three if $A$ is Dynkin type.

\end{minipage}
\end{center}
\end{abstract}

\vskip0.1in


\vskip0.1in


\footnote{ $^\star$Supported by the NSF of Shandong Province (Grant
No. Y2008A05)}

\footnote{ $^*$Corresponding author}

\footnote{  Email addresses: \ yinhongbo0218@126.com(H.Yin),
 \ shzhang@sdu.edu.cn(S.Zhang)}

\section{Introduction}

Representation dimension of Artin algebras was introduced by M.
Auslander in \cite{Auslander}, this concept gives a reasonable way
of measuring how far an Artin algebra $\Lambda$ is from being
representation-finite. In particular, M.Auslander has shown that an
Artin algebra is representation-finite if and only if its
representation dimension is at most 2.

\vskip0.2in

O.Iyama, in \cite{Iyama},  proved that the representation dimension
of an Artin algebra is always finite. Recently, Rouquier proved
in~\cite{Rouquier} that the representation dimension of Artin
algebras can be arbitrary large.

\vskip0.2in

An interesting relationship between the representation dimension and
the finitistic dimension conjecture has been shown by K. Igusa and
G. Todorov \cite{IT}, which is, if the representation dimension of
an algebra is at most three, then its finitistic dimension is
finite. Since then, many important algebras were proved to have
representation dimension at most three. Such as tilted algebras,
$m$-replicated algebras,  qusi-tilted algebras etc.,
see~\cite{APT}\cite{zl}\cite{Oppermann}for details.

\vskip0.2in

We follow this direction and investigate some special kinds of
triangular matrix algebras with small representation dimensions.

\vskip0.2in

Let $A$ be a finite dimensional hereditary algebra over an
algebraically closed field $k$, $T_2(A)=\left(\begin{array}{cc}A&0\\
A&A\end{array}\right)$ be the triangular matrix algebra and
$A^{(1)}=\left(\begin{array}{cc}A&0\\ DA&A\end{array}\right)$ be the
duplicated algebra of $A$ respectively.

\vskip0.2in

The following theorems are the main results of this paper.

\vskip0.2in

{\bf Theorem 1.}\ \ {\it Let $A$ be a finite dimensional hereditary
algebra over an algebraically closed field $k$. Then  ${\rm
rep.dim}\ T_2(A)\leq 3$ if $A$ is Dynkin type and $3\leq {\rm
rep.dim}\ T_2(A)\leq 4$ if $A$ is not Dynkin type.}

\vskip0.2in

{\bf Remark.}\  Theorem 1 improves the well known result about
representation dimension of $T_2(A)$.  According to \cite{FGR}, we
know that ${\rm rep.dim}\ T_2(A)\leq {\rm rep.dim}\ A+2$, which
implies that ${\rm rep.dim}\ T_2(A)\leq 5$ if $A$ is a finite
dimensional hereditary algebras over an algebraically closed field.

\vskip0.2in

Tilting theory of duplicated algebra $A^{(1)}$ has strong
relationship with cluster tilting theory induced in \cite{BM}, and
it has been widely investigated in \cite{ABSG, llz, z1, z2}. In this
paper, we mainly investigate the representation type and
representation dimension of endomorphism algebras of tilting modules
over duplicated algebra $A^{(1)}$.

\vskip0.2in

Let $T$ be a basic tilting A-$\module$ and $\ol{T}=T\oplus\ol{P}$ be
a tilting $A^{(1)}$-module, where $\ol{P}$ is the direct sum of all
non-isomorphic indecomposable projective-injective
$A^{(1)}$-modules.

\vskip0.2in

{\bf Theorem 2.}\ \ {\it Take the notation as above.  Then
$\End_{A^{(1)}}~\ol{T}$ is representation finite if and only if the
full subcategory $\{(X,Y,f)\ |\ X\in{\rm mod} A, \
Y\in\tau^{-1}\mathscr{F}(T_A)\cup{\rm add}\ A\}$ of ${\rm mod\
T_2(A)}$ is of finite type, where $\tau$ is the
Auslander-Reiten translation and $\mathscr{F}(T_A)$ is the
torsion-free class associated with $T$.}

\vskip0.2in

{\bf Theorem 3.}\ \ {\it Take the notation as above and assume that
$A$ is Dynkin type. Then ${\rm rep.dim\ End}_{A^{(1)}}\ {\ol T}\leq
3$.}

\vskip0.2in

{\bf Remark}\ We should mention that Theorem 1 can be obtained from
Theorem 3 by taking $T=DA$. We prove them differently, which seems
to be of independent interest.

\vskip 0.2in

This paper is arranged as follows. In Section 2, we collect
definitions and basic facts needed for our research. Section 3 is
devoted to the proof of Theorem 1, and in section 4, we prove
Theorem 2 and Theorem 3.

\vskip 0.2in

\section{Preliminaries}

\vskip0.2in

Let $\Lambda$ be a finite dimensional algebra over an
algebraically closed field $k$. We denote by mod $\Lambda$ the
category of all finitely generated right $\Lambda$-modules and by
ind $\Lambda$ the full subcategory of mod $\Lambda$ containing exactly
one representative of each isomorphism class of indecomposable
$\Lambda$-modules. We denote by ${\rm pd}\ X$ (resp. ${\rm id}\ X$)
the projective (resp. injective) dimension of an $\Lambda$-module
$X$ and by ${\rm gl.dim}\ \Lambda$ the global dimension of
$\Lambda$. Let $D={\rm Hom}_k(-,\ k)$ be the standard duality
between mod $\Lambda$ and mod $\Lambda^{{\rm op}}$, and
$\tau_\Lambda$ be the Auslander-Reiten translation of $\Lambda$. The
Auslander-Reiten quiver of $\Lambda$ is denoted by $\Gamma_\Lambda$.

\vskip0.2in

Let $M$ be a  $\Lambda$-module. We denote by ${\rm add}\ M$ the
subcategory of ${\rm mod}\ \Lambda$ whose objects are the direct
summands of finite direct sums of $M$. A module $M$ is called a
generator if all projective modules are in ${\rm add }\ M$ and is
called a cogenerator if all injective modules are in ${\rm  add}\
M$. We denote by ${\rm rep.dim}\ \Lambda$ the representation
dimension of $\Lambda$ which is defined by Auslander
in \cite{Auslander} as following.
$$
{\rm rep.dim}\ \Lambda={\rm min}\{ \ {\rm gl.dim}\
 {\rm End}_{\Lambda} M\ |\   M \ {\rm is \ a \ generator-cogenerator\
 for} \ \Lambda\}
$$

\vskip0.2in

The generator-cogenerator realizing the representation dimension is
called the Auslander generator. The following well known lemma is
the crucial tool to determine the upper bound of the representation
dimension.

\vskip0.2in

{\bf Lemma 2.1.}\ \cite{Auslander,EHIS,Xi}. {\it Let $A$ be an Artin
algebra. $M$ is a generator-cogenerator for $\Mod~A$. Then
$\gld~\End_A~M\leq n+2$ if and only if for each $A$-$\module$ $X$
there is an exact sequence
$$
(*)\ \ \
\xymatrix@C=0.5cm@R=0.15cm{0\ar[r]&M_n\ar[r]&\cdots\ar[r]&M_{0}\ar[r]&X\ar[r]&0}
$$
with all $M_i$  belongs to $\add~M$, such that the induced sequence
 $$\xymatrix@C=0.5cm@R=0.15cm{0\ar[r]&\Hom_A(M,M_n)\ar[r]&\cdots\ar[r]&\Hom_A(M,M_0)\ar[r]&\Hom_A(M,X)\ar[r]&0}$$
is exact.}

\vskip0.2in

{\bf Remark.}\ The exact sequence of $(*)$ in Lemma 2.1 is called an
${\rm add}\ M$-resolution of $X$.

\vskip0.2in

Let $\cal C$ be an additive Krull-Schmit Hom-finite $k$-category,
and $\cal X$ a full subcategory of $\cal C$.  We denote by ${\rm
ind}\ \cal X$ the subcategory consisting of indecomposable objects
of $\cal X$ and we say $\cal X$ is of finite type if ${\rm ind}\
\cal X$ is a finite set. Recall from \cite{AR}, a map $X'\rightarrow
A$ with $X'\in\cal X$ and $A\in\cal C$ is called a right $\cal
X$-approximation of $A$ if the induced map ${\rm
Hom}(X,X')\rightarrow {\rm Hom}(X,A)$ is an epimorphism for all
$X\in\cal X$. A map $f: A\rightarrow B$ in category $\cal C$ is
called right minimal, if for every $g: A\rightarrow A$ such that
$fg=f$, the map $g$ is an isomorphism. A right (left) approximation
that is also a right (left) minimal map is called a minimal right
(left) approximation of $T$. The subcategory $\cal X$ is called
contravariantly finite if any object in $\cal C$ admits a (minimal)
right $\cal X$-approximation. The notions of (minimal) left $\cal
X$-approximation and covariantly finite subcategory can be defined
dually.

\vskip 0.2in

A module $T\in \Mod\ \Lambda$ is called a tilting module if
the following conditions are satisfied:\\
{(1)} ${\rm pd}_\Lambda T \leq 1$;\\
{(2)} ${\rm Ext}_\Lambda^{1}(T,T)=0 $;\\
{(3)} There is  an exact sequence
 $0\longrightarrow \Lambda \longrightarrow T_{0}\longrightarrow T_{1}
\longrightarrow 0$ with $T_{i}\in {\rm add}\ T $ for $0\leq i\leq
1$.

\vskip 0.2in

Let $T$ be a tilting $\Lambda$-module and $B={\rm End}_\Lambda\ T$.
According to \cite{HR},  $(\mathscr{T}(T),\mathscr{F}(T))$ is the
torsion pair in $\Mod\ \Lambda$ generated by $T$, where
$\mathscr{T}(T)= T^{\bot}=\{\ X\in{\rm mod}\ \Lambda \ |\ {\rm
Ext}^1_{\Lambda}(T, X)=0\ \} = {\rm gen}\ T$ and
$\mathscr{F}(T)=\{X\in \Mod\ A\ |\ \ {\rm Hom}_\Lambda(T,X)=0
\}$, the corresponding torsion pair in $B$-mod is
$(\mathscr{X}(T),\mathscr{Y}(T))$, where $\mathscr{X}(T)=\{X\in
\Mod\ B\ |\ \ T\otimes_B X=0 \}$ and $\mathscr{Y}(T)=\{Y\in
\Mod\ B\ |\ \ {\rm Tor}^B_1(T,Y)=0 \}$.

\vskip 0.2in

A torsion pair is called splitting if every indecomposable
$\Lambda$-module  either belongs to the torsion class or belongs to
the torsion-free class. A tilting module $T$ is called splitting if
the corresponding torsion pair $(\mathscr{T}(T),\mathscr{F}(T))$ is
splitting, and $T$ is called separating if the corresponding torsion
pair $(\mathscr{X}(T),\mathscr{Y}(T))$ is splitting. Note that every
tilting module of hereditary algebras is splitting.

\vskip0.2in

{\bf Lemma 2.2.}\ {\it Let T be a tilting module of algebra
$\Lambda$. $B=\End_{\Lambda}~ T$}.

\vskip0.1in

(i)\ {\it $\Hom_\Lambda(T,-): \mathscr{T}(T)\rightarrow
\mathscr{Y}(T)$ and $-\otimes_{B}T: \mathscr{F}(T)\rightarrow
\mathscr{X}(T)$ are equivalent functors.}

\vskip0.1in

(ii)\  {\it $T$ is splitting if and only if ${\rm id}\ X =1$ for
every $X \in\mathscr{F}(T)$}

\vskip0.1in

(iii) {\it $T$ is separating if and only if ${\rm pd} \ Y=1$ for
every $Y\in\mathscr{X}(T)$.}

\vskip0.2in

 Let ${\cal T}_{\Lambda}$ be the set of
all basic tilting $\Lambda$-modules up to isomorphism. Recall from
\cite{HU}, the tilting quiver $\mathscr{K}(\Lambda)$ of $\Lambda$ is
defined as the following. The vertices of $\mathscr{K}(\Lambda)$ are
the elements of ${\cal T}_{\Lambda}$. There is an arrow
$T'\rightarrow T$ in $\mathscr{K}(\Lambda)$ if and only if
$T'=M\oplus X$ and $T=M\oplus Y$ with $X$ and $Y$ indecomposable
such that there is a short exact sequence $0\rightarrow X\stackrel
{f}\longrightarrow E\stackrel {g}\longrightarrow Y\rightarrow 0$
such that $f$ is a minimal left ${\rm add}\ M$-approximation of $X$
and that $g$ is a minimal right ${\rm add}\ M$-approximation of $Y$.

\vskip0.2in

We recall the definition of triangular matrix algebra from
\cite{ARS}. Let A and B be finite dimensional algebras over $k$ and
$_AM_B$ be an $A$-$B$-bimodule. $\Lambda=\left(\begin{array}{cc}B&0\\
M&A\end{array}\right)$ is called a triangular matrix algebra, its
elements are $\left(\begin{array}{cc}b&0\\ m&a\end{array}\right)$
where $b\in B, a\in A, m\in M$, and its addition and multiplication
are given by the usual matrix operation.

\vskip0.2in

{\bf Remark.}\ There are two special kinds of
triangular matrix algebras. One is $T_2(A)=\left(\begin{array}{cc}A&0\\
A&A\end{array}\right)$ and the other is $A^{(1)}=\left(\begin{array}{cc}A&0\\
DA&A\end{array}\right)$which is also called the duplicated algebra
of $A$

\vskip0.2in

It is well known
that the module category of $\Lambda=\left(\begin{array}{cc}B&0\\
M&A\end{array}\right)$ is equivalent to the category
$\mathrm{rep}(_AM_B)$, called the representation of the bimodule
$_AM_B$, see \cite[AppendixA 2.7]{ASS}. The objects of
$\mathrm{rep}(_AM_B)$ are triples $(X, Y, f)$, where $X$ is an
$A$-$\module$, $Y$ is a $B$-$\module$ and $f: X\otimes_A
M\rightarrow Y$ is a $B$-$\module$ morphism. The morphism between
$(X_1, Y_1, f_1)$ and $(X_2, Y_2, f_2)$ is a pair $(x, y)$ makes the
following diagram commutative.

$$
\xymatrix{
          X_1\otimes_A M\ar[r]^{x\otimes_A M}\ar[d]_{f_1}&X_2\otimes_A M\ar[d]_{f_2}\\
          Y_1\ar[r]^{y}&Y_2
          }
$$
Using the adjoint isomorphism between $-\otimes M$ and $\Hom(M,-)$,
the category $\mathrm{rep}(_AM_B)$ can also be described as follows.
Its objects are triples $(X, Y, f)$, where X is an $A$-$\module$, Y
is a $B$-$\module$ and $f: X \rightarrow \Hom_B(M,Y)$ is an
$A$-$\module$ morphism. The morphism between $(X_1, Y_1, f_1)$ and
$(X_2, Y_2, f_2)$ is a pair $(x, y)$ makes the following diagram
commutative.

$$
\xymatrix{
          X_1\ar[r]^{x}\ar[d]_{f_1}&X_2\ar[d]_{f_2}\\
          \Hom_B(M,Y_1)\ar[r]^{\scriptscriptstyle\Hom(M,y)}&{\Hom_B(M,Y_2)}
          }
$$
In the following, we will freely use the two descriptions as the modules
of the triangular matrix algebras.

\vskip0.2in

The indecomposable projective $\Lambda$-$\modules$ are isomorphic to
objects of the form $(0,Q,0)$ where $Q$ is an indecomposable
projective $B$-$\module$ and $(P, P\otimes_A M, {\rm id})$ where $P$
is an indecomposable $A$-$\module$, ${\rm id}$ is the identity map.
Dually, the indecomposable injective $\Lambda$-$\modules$ are
isomorphic to objects of the form $(I,0,0)$ where $I$ is an
indecomposable injective $A$-$\module$ and $(\Hom_B(M,J),J,{\rm
id})$ where $J$ is an indecomposable injective $B$-$\module$. See
\cite[III, Proposition 2.5]{ARS} for details.

\vskip0.2in

Let $\cal C$ be an additive Krull-Schmit Hom-finite $k$-category. We
define $\Mor~\cal C$ to be the morphism category of $\cal C$ whose
objects are triples $(X, Y, f)$ where $X, Y\in{\cal C},
f\in\Hom_{\cal C}(X,Y)$ and the morphism between $(X_1, Y_1, f_1)$
and $(X_2, Y_2, f_2)$ is a pair $(x, y)$ makes the following diagram
commutative.

$$
\xymatrix{
          X_1\ar[r]^{x}\ar[d]_{f_1}&X_2\ar[d]_{f_2}\\
          Y_1\ar[r]^{y}&Y_2
          }
$$
We also use $\xymatrix@C=0.5cm{X\ar[r]^f&Y}$ to denote the objects
of $\Mor~\cal C$. It is easy to see that $\Mor~\cal C$ is also a
Krull-Schmidt category. In particular, $\Mor({\rm mod}\ A)$ is
equivalent to ${\rm mod}\ T_2(A)$.

\vskip0.2in

Let ${\cal X}_1$ and ${\cal X}_2$ be two full subcategories of $\cal
C$. We denote by ${\cal X}_1\cup{\cal X}_2$ the full subcategory of
$\cal C$ consisting of the objects which either belong to ${\cal
X}_1$ or belong to ${\cal X}_2$ and ${\cal X}_1\setminus{\cal X}_2$
the full subcategory of $\cal C$ consisting of the objects belonging
to ${\cal X}_1$ and not belonging to ${\cal X}_2$. If ${\cal C}={\rm
mod}\ \Lambda$ is a module category, then we denote by
$\tau\mc_1$(resp. $\Tau\mc_1$) the full subcategory of $\mc$ whose
objects are obtained from $\mc_1$ by the once action of $\tau$(resp.
$\Tau$).

\vskip0.2in

Let $\Lambda$ be an Artin algebra of finite representation type and
$M_1, M_2,\ldots, M_n$ be a complete set of non-isomorphic
indecomposable $\Lambda$-modules. According to Auslander in
\cite{ARS}, $M=M_1\oplus\cdots\oplus M_n$ is called an additive
generator of ${\rm mod}\ \Lambda$ and $\Gamma_M=\End_A~(M)^{{\rm
op}}$ is said to be the Auslander algebra of $\Lambda$. It is well
known that $\Gamma_M$ and $T_2(\Lambda)$ have the same
representation type.

\vskip0.2in

Let $\mc$ be an additive Krull-Schmidt $\Hom$-finite k-category of
finite type and $M$ be an object of $\mc$. $M$ is said to be an
additive generator of $\mc$ if every indecomposable object of $\mc$
is a direct summand of $M$, then $\Gamma_M$=${{\rm End}_{\mc}\
(M)}^{{\rm op}}$ is said to be the Auslaner algebra of $\mc$.

\vskip0.2in

{\bf Remark.}\ \ If $\mc={\rm mod}\ \Lambda$ for some
representation-finite algebra $\Lambda$, then the Auslander algebra
of $\mc$ is the same as the original definition of Auslander algebra
of $\Lambda$.

\vskip0.2in

The following proposition is similar with Proposition 5.8 in
\cite[p.215]{ARS}.

\vskip0.2in

{\bf Proposition 2.3.}\ {\it Let $\mc$ be an additive Krull-Schmidt
$\Hom$-finite k-category of finite type and M an additive generator
of $\mc$. Then the Auslander algebra $\Gamma_M$ of $\mc$ is of
finite representation type if and only if the category
$\mathrm{Mor}~\mc$ is of finite type.}

\vskip0.2in

Throughout this paper, the notations will be fixed as above. We
refer to \cite{ASS, ARS} for the other concepts of representation
theory of Artin algebras.

\vskip0.2in

\section{Representation dimension of triangular matrix algebras}

\vskip0.2in

In this section, we assume that $A$ is a finite dimensional
hereditary algebra over an algebraicall closed filed $k$. We will
give a bound of the representation dimension of $T_2(A)$  by using
Proposition 2.3, and then prove Theorem 1.

\vskip0.2in

{\bf Theorem 3.1.}\ {\it  Let $ A$ be a hereditary algebra of Dynkin
type, then ${\rm rep.dim}\ T_2( A)\leq 3$.}

\vskip0.1in

{\bf Proof}\ Note that $A$ is an additive generator of ${\rm add}\
A$ and that $\End_A A$=$A^{\rm op}$ is of finite representation
type. Then $\Mor~({\rm add}\ A)$ is of finite type by Proposition
2.3. By the same argument we know that
$\Mor~(\add\ DA)$ is also of finite type.

Let $(P_1, P_2, f)$ and $(I_1, I_2, g)$ be additive generators of
$\Mor~({\rm add}\ A)$ and $\Mor~\add\ DA$ respectively and let
$M=(P_1, P_2, f)\oplus(I_1, I_2, g)$. Then $M$ is a
generator-cogenerator for ${\rm mod}\ T_2( A)$.

We claim that ${\rm gl.dim\ End}_{T_2( A)} M\leq 3$.

In fact, let $(X, Y, h)$ be an indecomposable $T_2( A)$-module. We
may assume that $(X, Y, h)$ does not belong to $\add~M$.

\vskip0.1in

{\bf Case I.}\ Assume that $X$ and $Y$ have no non-zero injective
direct summand.  Then there is a minimal right
$\add~M$-approximation $(\alpha_1, \alpha_2): (P'_1, P'_2,
f')\rightarrow (X, Y, h)$ of $(X, Y, h)$ with $(P'_1, P'_2, f')$
belongs to $\Mor~\mathscr{P}( A)$ since there is no non-zero
homomorphism from $(I_1, I_2, g)$ to $(X, Y, h)$. Then we have the
following commutative diagram.
\begin{equation}
\begin{array}{c} \xymatrix{
          0\ar[r]&P''_1\ar[r]^{\beta_1}\ar[d]_{f''}&P'_1\ar[r]^{\alpha_1}\ar[d]_{f'}&X\ar[d]_{h}\ar[r]&0\\
          0\ar[r]&P''_2\ar[r]^{\beta_2}&P'_2\ar[r]^{\alpha_2}&Y\ar[r]&0
          }
\end{array}
\end{equation}
with exact rows. This is a short exact sequence of $T_2(
A)$-modules, and it follows that $P''_1$ and $P''_2$ are projective
$A$ modules since $A$ is hereditary and $P'_1, P'_2$ are projective.
In particular, $(P'_1, P''_2, f'')$ belongs to $\add~M$ and the
diagram of $(3.1)$ is an $\add~M$-resolution of $(X, Y, h)$.

{\bf Case II.}\ $X$ or $Y$ have non-zero injective direct summand,
that is, $(X, Y, h)$ is of the form $h=\left(\begin{array}{cc}
h_1&0\\h_2&h_3\end{array}\right):\xymatrix{X_1\oplus I'_1\ar[r]&Y_1
\oplus I'_2}$, where $X_1$ and $Y_1$ both have no non-zero injective
direct summand,  $I_1'$ and $I_2'$ are injective $A$-modules (may be
zero).

Let $i_1=\left(\begin{array}{c} 0\\1\end{array}\right):
I_1'\rightarrow X_1\oplus I_1'$ and $i_2=\left(\begin{array}{c}
0\\1\end{array}\right): I_2'\rightarrow Y_1\oplus I_2'$. It is easy
to see that
$$
\xymatrix{I'_1\ar[r]^-{i_1}\ar[d]_{h_3}&X_1\oplus
I'_1\ar[d]^{\scriptsize \left(\begin{array}{cc}
h_1&0\\h_2&h_3\end{array}\right)}\\I'_2\ar[r]^-{i_2}&Y_1\oplus I'_2}
$$
is a minimal right $\add(I_1, I_2, g)$-approximation of $(X, Y, h)$.

Let
$$\xymatrix{P'_1\ar[r]^{\alpha_1}\ar[d]_{f'}&X\ar[d]_h\ar[r]&0\\P'_2\ar[r]^{\alpha_2}&Y\ar[r]&0}$$
be a minimal right $\add(P_1, P_2, f)$-approximation of $(X, Y, h)$.
It is epimorphism since all the indecomposable projective
$T_2(A)$-modules belong to $\add~(P_1, P_2, f)$. Then $(P'_1, P'_2,
f')\oplus (I'_1, I'_2, h_3)$ is a right $\add~M$-approximation of
$(X, Y, h)$.

We have the following commutative diagram.
\begin{equation}\label{1}
\begin{array}{c}\xymatrix{0\ar[r]&\Ker~\pi_1\ar[r]\ar[d]_\varphi&P'_1\oplus I'_1\ar[r]^-{\pi_1}\ar[d]_{\scriptsize
\left(\begin{array}{cc}
f'&0\\0&h_3\end{array}\right)}&X\ar[d]_{h}\ar[r]&0\\0\ar[r]&\Ker~\pi_2\ar[r]&P'_2\oplus
I'_2 \ar[r]^-{\pi_2}&Y\ar[r]&0}\end{array}\end{equation} with
$\pi_1=(\alpha_1, i_1)$ and $\pi_2=(\alpha_2, i_2)$. Now we need to
determine $\Ker~\pi_1$ and $\Ker~\pi_2$.

Consider the pull-back of $(\alpha_1, i_1)$:
$$
\xymatrix{D\ar[r]^{d_1}\ar[d]_{d_2}&I'_1\ar[d]_{i_1}\\P'_1\ar[r]^{\alpha_1}&X}
$$
which implies that $\Ker~\pi_1\simeq D$. Note that a pull-back
diagram is also a push-out diagram if and only if
$(\begin{array}{cc}\alpha_1&i_1\end{array})$ is epimorphism see
\cite[exercise 6.7]{HS}. Hence, the above diagram is also a push-out
of $(d_1, d_1)$. In particular, we have that $d_2$ is monomorphism
because $i_1$ is and $\Ker~\pi_1$ is projective. Then $\Ker~\pi_2$
is also projective by the same argument. Hence $(\Ker~\pi_1,
\Ker~\pi_2, \varphi)$ belongs to $\add~M$ and (\ref{1}) is an
$\add~M$-resolution of $(X, Y, h)$.

Summary the above discussions, we know that ${\rm gl.dim\ End}_{T_2(
A)} M\leq 3$, which forces that ${\rm rep.dim}\ T_2( A)\leq 3$. The
proof is completed. $\hfill\Box$

\vskip0.2in

{\bf Theorem 3.2.}\ {\it  Let $ A$ be a hereditary algebra of
Euclidean or wild type. Then $3\leq {\rm rep.dim}\ T_2( A)\leq 4$}

\vskip0.1in

{\bf Proof}\ The first inequality is obvious since $T_2(A)$ is
representation infinite when $A$ is not Dynkin type.

For the second inequality, we choose $M=T_2( A)\oplus DT_2( A)$ as a
generator-cogenerator for ${\rm mod}\ T_2( A)$. Note that $\gld~
T_2( A)\leq\gld~ A+1=2$, see for example \cite [Proposition 2.6, p
78]{ARS} for details.

Let $(X, Y, h)$ be a $T_2( A)$-module. We may assume that $(X, Y,
h)$ does not belong to $\add~M$.

If both $X$ and $Y$ have no non-zero injective direct summand, then
a projective resolution of $(X, Y, h)$  also is its
$\add~M$-resolution which obviously has length at most 2.

If $X$ or $Y$ have non-zero injective direct summand, then $(X, Y,
h)$ is of the form $\left(\begin{array}{cc}
h_1&0\\h_2&h_3\end{array}\right):\xymatrix{X_1\oplus
I'_1\ar[r]&Y_1\oplus I'_2}$, where $X_1$ and $Y_1$ both have no
non-zero injective direct summand,  $I_1'$ and $I_2'$ are injective
$A$ modules (may be zero).

Let $i_1=\left(\begin{array}{c} 0\\1\end{array}\right): I_1'\rightarrow
X_1\oplus I_1'$ and $i_3=\left(\begin{array}{c}
0\\h_3\end{array}\right): I_1'\rightarrow Y_1\oplus I_2'$.

Then
$$\xymatrix{I'_1\ar[r]^-{i_1}\ar[d]_{1}&X_1\oplus I'_1\ar[d]^{\scriptsize \left(\begin{array}{cc}
h_1&0\\h_2&h_3\end{array}\right)}\\I'_1\ar[r]^-{i_3}&Y_1\oplus
I'_2}$$ is a minimal right $\add~DT_2(A)$-approximation of $(X, Y,
h)$, and let
$$\xymatrix{P'_1\ar[r]^{\alpha_1}\ar[d]_{f'}&X\ar[d]_h\ar[r]&0\\P'_2\ar[r]^{\alpha_2}&Y\ar[r]&0}$$
is a minimal right $\add~T_2(A)$-approximation of $(X, Y, h)$. Then
we have a commutative diagram similar to (\ref{1})
\begin{equation}\label{2}
\begin{array}{c}\xymatrix{0\ar[r]&\Ker~\pi_1\ar[r]\ar[d]_\varphi&P'_1\oplus I'_1\ar[r]^-{\pi_1}\ar[d]_{\scriptsize
\left(\begin{array}{cc}
f'&0\\0&1\end{array}\right)}&X\ar[d]_{h}\ar[r]&0\\0\ar[r]&\Ker~\pi_2\ar[r]&P'_2\oplus
I'_1 \ar[r]^-{\pi_2}&Y\ar[r]&0\ .}\end{array}\end{equation}

By the same argument as in Theorem 3.1, we know $\Ker~\pi_1$ is
projective. Let $P''_2$ be the projective cover of $\Ker~\pi_2$.
Then $(\Ker~\pi_1, \Ker~\pi_1, 1)\oplus(0, P''_2, 0)$ is a right
$\add~T_2(A)$-approximation of $(\Ker~\pi_1, \Ker~\pi_2, \varphi)$,
and we have the following commutative diagram with exact rows:

\begin{equation}\label{3}
 \begin{array}{c}
  \xymatrix{0\ar[r]&0\ar[r]\ar[d]&\Ker~\pi_1\oplus0\ar[r]^-{1}\ar[d]_{\scriptsize\left(\begin{array}{cc}1&0\\
  0&0\end{array}\right)}&\Ker~\pi_1\ar[r]\ar[d]_{\varphi}&0\\
           0\ar[r]&P_3\ar[r]&\Ker~\pi_1\oplus P''_2\ar[r]^-{\scriptsize(\varphi\,p)}&\Ker~\pi_2\ar[r]&0\ ,}
 \end{array}
\end{equation}
where $P_3$ is a projective $A$-module. According to diagrams
(\ref{2}) and (\ref{3}), we get a right $\add~M$-resolution of $(X,
Y, h)$ of length at most 2, hence ${\rm rep.dim}\ T_2( A)\leq 4$.
This completes the proof.      $\hfill\Box$

\vskip0.2in

\section{Endomorphism algebras of tilting modules of duplicated algebras}

Let $A$ be a hereditary algebra and $T$ be a tilting right
$A$-$\module$. Let $B=\End_AT$. Then $_BT_A$ is a $B$-$A$-bimodule.
Let $A^{(1)}$ be the duplicated algebra of A and $\overline{P}$ be
the direct sum of all non-isomorphic indecomposable
projective-injective $A^{(1)}$-modules. Then
$\overline{T}=T\oplus\overline{P}$ is a tilting right
$A^{(1)}$-module. We will prove Theorem 2 and Theorem 3 in this
section.

\vskip0.2in

Note that $\overline{P}$ and $T$ regarded as $A^{(1)}$-modules can
be written as ${\overline P}=(A, DA, id)$ and $T=(0,T,0)$. Then we
have follows.
$$
\Hom_{A^{(1)}}((0,T,0),\overline{P})=\Hom_A(T,DA)=DT
$$
$$
\End_{A^{(1)}}~\ol{T}=\left(\begin{array}{cc}\End_AT&0\\\Hom_{A^{(1)}}((0,T,0),
\overline{P})&\End_{A^{(1)}}~\overline{P}\end{array}\right)
=\left(\begin{array}{cc}B&0\\DT&A\end{array}\right).
$$

\vskip0.2in

{\bf Theorem 4.1.}\ {\it  Take the notations as above. Then
$\End~\ol{T}$ is representation finite if and only if the full
subcategory $\{(X,Y,f)\ |\ X\in{\rm mod}\ A,
Y\in\tau^{-1}\mf(T_A)\cup{\rm add}\ A \}$ of ${\rm mod \
T_2(A)}$ is of finite type.}

\vskip0.2in

In order to prove the theorem, we need following lemmas.

\vskip0.2in

{\bf Lemma 4.2.}\ {\it Let $A$, $T$ and $B$ be as above. Then $D T$
is a separating convex tilting right $B$-module.}

\vskip0.1in

{\bf Proof}\   It is shown in \cite{HR} that ${\rm ind \  add}(D T)$
is a complete slice in ${\rm mod}\ B$. Hence $D T$ is a convex
tilting right $B$-module, and $D T$ is a splitting tilting left
$A$-module since $A$ is hereditary. Then it follows that $D T$ is
also a separating tilting right $B$-module.      $\hfill\Box$

\vskip0.2in

{\bf Lemma 4.3.}\ {\it Assume $A$, $T$ and $B$ be as above. Let
$(\mathscr{T}(T), \mathscr{F}(T))$ and $(\mathscr{X}(DT),
\mathscr{Y}(DT))$ be the torsion pairs in ${\rm mod}\ A$
corresponding to $T$ and $DT$ respectively, and let
$(\mathscr{T}(DT), \mathscr{F}(DT))$ and $(\mathscr{X}(T),
\mathscr{Y}(T))$ be the torsion pairs in ${\rm mod}\ B$
corresponding to $DT$ and $T$ respectively. Then we have the
following.}

\vskip0.1in

(i)\ $\mathscr{T}(DT)=\mathscr{X}(T)\cup {\rm add}\ DT$, \
$\mathscr{F}(DT)=\mathscr{Y}(T)\setminus {\rm add}\ DT$;

\vskip0.1in

(ii)\  \ $\mathscr{Y}(DT)=\tau^{-1}\mathscr{F}(T)\cup {\rm add}\ A$,
\ $\mathscr{X}(DT)=\tau^{-1}(\mathscr{T}(T)\setminus {\rm add}\ DA)$.

\vskip0.1in

{\bf Proof}\ (i) According to Lemma 4.2, $DT$ is a separating convex
tilting right $B$-module. We have the following.
$$
{\rm ind}\ \mathscr{T}(DT)=\{M\in {\rm ind}\ B\ |\ {\rm \Hom}_B(DT,
M)\neq 0\},
$$
$$
{\rm ind}\ \mathscr{F}(DT)=\{M\in {\rm ind}\ B\ |\ {\rm \Hom}_B(DT,
M)=0\}.
$$
Let $M$ be an indecomposable $B$-module.  Then either
$M\in\mathscr{Y}(T)$ or $M\in\mathscr{X} (T)$ since $T$ is
splitting.

If $M\in\mathscr{Y}(T)$, then $M=\Hom_A(T, N)$ for some $N\in\ind\
\mathscr{T}(T)$,
$$
\Hom_B(DT, M) =\Hom_B(\Hom_A(T, DA), \Hom_A(T, N))=\Hom_A(DA, N).
$$
In this situation, $M\in\mathscr{T}(DT)$ if and only if
$N\in\add~DA$, which equivalent to $M\in \add~DT$.

If $M\in\mathscr{X}(T)$, then $M=\Ext(T, N)$ for some $N\in\ind \
\mathscr{F}(T)$. Applying $\Hom_A(T,-)$ to an injective resolution
$\xymatrix{0\ar[r]&N\ar[r]&I_1\ar[r]&I_2\ar[r]&0}$ of $N$, we get
the following exact sequence
$$
\xymatrix@C=0.5cm{0\ar[r]&\Hom_A(T, N)\ar[r]&\Hom_A(T,
I_1)\ar[r]&\Hom_A(T, I_2)\ar[r]&\Ext(T, N)\ar[r]&0}.
$$
Therefore, $\Hom_B(\Hom_A(T, I_2), \Ext(T, N))\neq 0$.

On the other hand, by using $DT=\Hom_A(T, DA)$ we have
$$
\Hom_B(DT, M)=\Hom_B(\Hom_A(T, DA), \Ext(T, N))\neq 0,
$$
hence $\mathscr{T}(DT)=\mathscr{X}(T)\cup\add~ DT$. Since $DT$ is
separating, we have that
$\mathscr{F}(\dt)=\mathscr{Y}(T)\setminus\add~DT$.

\vskip0.1in

(ii)\  We have $\mathscr{Y}(\dt)=\{\Hom_B(DT, M) \ | \
M\in\mathscr{T}(DT)\}$, and by (i), $\mathscr{Y}(\dt)=\{\Hom_B(DT,
M) \ | \ M\in\mathscr{X}(T)\cup {\rm add}\ DT\}$.

If $M$ is an indecomposable $B$-module which belongs to $\add~\dt$,
then $\Hom_B(\dt, M)\in\add~ A$.

If $M\in{\rm ind}\ \mathscr{X}(T)$, then $M=\Ext(T, N)$ for some
$N\in{\rm ind}\ \mathscr{F}(T)$. Then we have
\begin{eqnarray*}
\Hom_B(\dt, M)&=&\Hom_B(\Hom_A(T, DA), \Ext(T, N))\\
              &=&\Ext(DA, N)\\
              &=&D\Hom_A(\Tau N, DA)\\
              &=&\Tau N
\end{eqnarray*}
By using the equivalence of $\mathscr{X}(T)$ and $\mathscr{F}(T)$,
we have $\{\Hom_B(\dt, M)\ |\
M\in\mathscr{X}(T)\}=\Tau\mathscr{F}(T)$. Hence,
$\mathscr{Y}(\dt)=\Tau\mathscr{F}(T)\cup \add~ A$.

\vskip0.1in

Finally, we determine  $\mathscr{X}(\dt)$. According to Lemma 2.2
and by using (i), we know that $\mathscr{X}(\dt)=\{
\mathrm{Ext}_B^1(\dt, M)\ |\ M\in\mathscr{F}(DT)\}=\{
\mathrm{Ext}_B^1(\dt, M)\ |\ M\in\mathscr{Y}(T)\backslash{\rm add}\
DT \}$.

Let $M$ be an indecomposable $A$-module in $\mathscr{X}(\dt)$. Then
$M=\Hom_A(T, N)$ for some  $N\in\ind \ \mathscr{T}(T)$ with
$N\notin\add~ DA$.

Then we have
\begin{eqnarray*}
\mathrm{Ext}_B^1(\dt, M)&=&\mathrm{Ext}_B^1(\dt, \Hom_A(T, N))\\
               &=&\Ext(\dt\otimes_B{T}, N)\\
               &=&\Ext(DA, N)\\
               &=&D\Hom_A(\Tau N, DA)\\
               &=&\Hom_A(A, \Tau N)\\
               &=&\Tau N
\end{eqnarray*}
Hence, $\mathscr{X}(\dt)$=$\Tau(\mathscr{T}(T)\setminus\add~ DA)$.
The proof is completed.              $\hfill\Box$

\vskip0.2in

{\bf Proof of Theorem 4.1.}\  \  If A is not Dynkin type, then both
$\{(X,Y,f)|X\in \Mod~A,\ Y\in\tau^{-1}\mathscr{F}(T_A)\cup{\rm add}\
A\}$ and ${\rm End}_{A^{(1)}}~\ol T$ are of infinite type. Hence
without loss of generality, we can assume that A is Dynkin type.

Now, let $(X, Y, g)$ be any indecomposable ${\rm End}_{A^{(1)}}~\ol
T$-module. According to Lemma 4.2, $DT$ is a separating convex
tilting $B$-module, hence $Y$ can be written as $Y=M\oplus N$ where
$M\in\mathscr{T}(\dt)$ and $N\in\mathscr{F}(\dt)$. By Lemma 4.3 (i)
$\Hom_B(\dt, N)$=0, and there are two kinds of indecomposable ${\rm
End}_{A^{(1)}}~\ol T$-modules.

(1) \ $N\neq 0$. Then $(0, N, 0)$ is a direct summand of $(X,
M\oplus N, g)$, it forces that $(X, M\oplus N, g)=(0, N, 0)$ and $N$
is indecomposable. Note that the number of indecomposable ${\rm
End}_{A^{(1)}}~\ol T$ modules of this kind is finite since A is
Dynkin type.

(2)  $N=0$. Then $(X, M\oplus N, g)$=$(X, M, g)$. Hence ${\rm
End}_{A^{(1)}}~\ol T$ is representation finite if and only if there
are finite number of indecomposable modules of the second type. We
write the full subcategory of the second kind of modules by
$\Mod_2~\End_{A^{(1)}}~\ol T$. We define a Functor $F:(X, M,
g)\rightarrow (X, \Hom_B(\dt, M), g')$. By Lemma 4.3 (ii) we know
that $F$ is an equivalence between the $\Mod_2~\End_{A^{(1)}}~\ol T$
and $\{(X, Y, f)\in\Mod~T_2(A)|\ X\in\Mod~A,\ Y\in\add~
A\cup\Tau\mf(T)\}$. This completes the proof of the theorem.
$\hfill\Box$

\vskip0.2in

{\bf Example.}\  Let $A$ be the path algebra of the quiver
$$
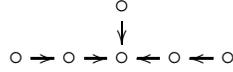
\xymatrix@C=0.3cm @R=0.3cm{&&\circ\ar[d]&&\\ \circ\ar[r]&\circ\ar[r]&\circ&\circ\ar[l]&\circ\ar[l]}$$
the Aulander-Retein quiver of $A$ is as follows.

$$
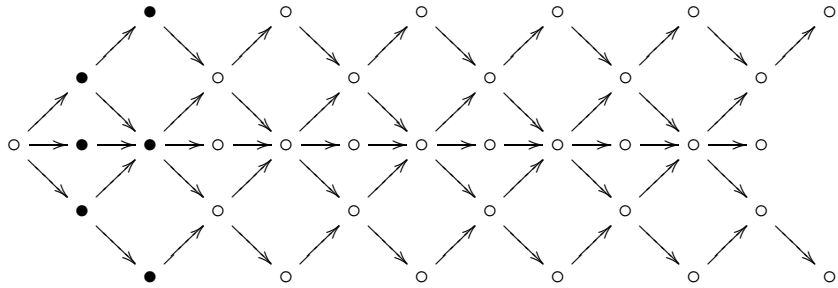
\xymatrix@C=0.5cm@R=0.5cm{&&\bullet\ar[rd]&&\circ\ar[rd]&&\circ\ar[rd]&&\circ\ar[rd]&&\circ\ar[rd]&&\circ\\
            &\bullet\ar[ru]\ar[rd]&&\circ\ar[ru]\ar[rd]&&\circ\ar[ru]\ar[rd]&&\circ
            \ar[ru]\ar[rd]&&\circ\ar[ru]\ar[rd]&&\circ\ar[ru]&\\
            \circ\ar[ru]\ar[rd]\ar[r]&\bullet\ar[r]&\bullet\ar[ru]\ar[rd]\ar[r]&\circ\ar[r]&\circ
            \ar[ru]\ar[rd]\ar[r]&\circ\ar[r]&\circ\ar[ru]\ar[rd]\ar[r]&\circ\ar[r]&\circ\ar[ru]\ar[rd]
            \ar[r]&\circ\ar[r]&\circ\ar[ru]\ar[rd]\ar[r]&\circ&\\
            &\bullet\ar[ru]\ar[rd]&&\circ\ar[ru]\ar[rd]&&\circ\ar[ru]\ar[rd]&&\circ\ar[ru]
            \ar[rd]&&\circ\ar[ru]\ar[rd]&&\circ\ar[rd]&\\  &&\bullet\ar[ru]&&\circ\ar[ru]&&
            \circ\ar[ru]&&\circ\ar[ru]&&\circ\ar[ru]&&\circ}$$
Consider the APR-tilting module $T$ of $A$ whose indecomposable
direct summand is denoted by $\bullet$ in the Auslander-Reiten
quiver. $\mathscr{F}(T)$ has one indecomposable module which is the
simple projective $A$-module. It is easy to see that
$\Tau\mf(T)\cup\add A$ is of finite type and its Auslander algebra
$\Gamma$ is given by the quiver
$$
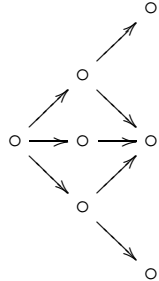
\xymatrix@C=0.5cm@R=0.5cm{&&\circ\\
                            &\circ\ar[ru]\ar[rd]\\
                            \circ\ar[ru]\ar[r]\ar[rd]&\circ\ar[r]&\circ\\
                            &\circ\ar[rd]\ar[ru]\\
                            &&\circ}$$
with the relation that the sum of the three roads with length two in
the middle mesh equals zero. Note that $\Gamma$ is representation
infinite, hence by Theorem 4.1 and Proposition 2.3,
$\End_{A^{(1)}}~\ol T$ is also representation infinite.

\vskip0.2in

{\bf Corollary 4.4.}\ {\it  Let $A$ be a finite dimensional
hereditary algebra over an algebraically closed field $k$ and let
$T_1$ and $T_2$ be two (basic) tilting module of $A$. Assume that
there is a path from $T_1$ to $T_2$ in the tilting quiver
$\mathscr{K}(A)$, if  $\End_{A^{(1)}}~\ol T_1$ is representation
infinite, then $\End_{A^{(1)}}~\ol T_2$ is also representation
infinite.}

\vskip0.1in

{\bf Proof}\  It follow from the fact that
$\mathscr{T}(T_1)\subseteq\mathscr{T}(T_2)$ since there is a path
from $T_1$ to $T_2$.      $\hfill\Box$

\vskip0.2in

{\bf Remark.}\  The converse of Corollary 4.4 is not true, that is,
there exists tilting modules $T_1$ and $T_2$ with a path  between
them in the tilting quiver such that $\End_{A^{(1)}}~\ol{T}_1$ is
representation finite, but $\End_{A^{(1)}}~\ol{T}_2$ is
representation infinite. See the example above, $A$ is
representation finite while the unique APR-tilting module is
representation infinite and there is an arrow from $A$ to the
APR-module.

\vskip0.2in

In the rest part of this section, we investigate the representation
dimension of Endomorphism algebras of tilting modules over
duplicated algebras. The following lemma is useful in our research.

\vskip0.2in

{\bf Lemma 4.5.}\ {\it  Let $A$ be a finite dimensional hereditary
algebra over an algebraically closed field $k$ and $T$ be a tilting
$A$-module. Assume that $f: T_1\rightarrow X$ is a right
$\add~ T$-approximation of X, then $\Ker f\in\add~ T$.}

\vskip0.1in

{\bf Proof}\  $\Hom_A(T,X)=\Hom_A(T, tX)$, hence $f$ is in fact the
composition of $\xymatrix@C=0.5cm{T_1\ar[r]^-{f}&tX\ar[r]&X}$. Note
$\xymatrix@C=0.5cm{T_1\ar[r]^{f}&tX}$ is the $\add~ T$-approximation
of $tX$. Because $tX\in\mathscr{T}(T)$, $tX$ is generated by $T$. So,
$\xymatrix@C=0.5cm{T_1\ar[r]^f&tX\ar[r]&0}$ is epimorphism. Applying
$\Hom_A(T,-)$ to the exact sequence
\begin{equation}\label{5}
\xymatrix@C=0.8cm{0\ar[r]&\Ker f\ar[r]&T_1\ar[r]&tX\ar[r]&0}
\end{equation}
We get
$$\xymatrix@C=0.6cm@R=0.15cm{\Hom_A(T,T_1)\ar[r]^-{\scriptsize (T,f)}&
\Hom_A(T,tX)\ar[r]&\Ext(T,\Ker f)\ar[r]&\Ext(T,T_1)=0}$$
Because $\Hom_A(T,f)$ is epimorphism, we have $\Ext(T,\Ker f)=0$.
So, $\Ker f\in \mathscr{T}(T)$.

Let $U$ be any module in $\mathscr{T}(T)$. Applying $\Hom(-,U)$ to (\ref{5})
we get
$$\xymatrix{0=\Ext(T_1, U)\ar[r]&\Ext(\Ker f, U)\ar[r]&\mathrm{Ext}^2_A(tX, U)=0}$$
So, $\Ext(\Ker f, U)=0$, i.e. $\Ker f$ is Ext-projective in
$\mathscr{T}(T)$. We know that $\Ker f$ belongs to $\add~ T$. This completes
the proof.      $\hfill\Box$

\vskip0.2in

The following lemma is taken from \cite[Proposition 2.2]{APT} which
will be used later.

\vskip0.2in

{\bf Lemma 4.6.}\ {\it Let A be an Artin algebra. $M=T\oplus N$ is
an A-module. X is generated by M and
$\xymatrix@C=0.5cm{0\ar[r]&K\ar[r]&M_0\ar[r]&X\ar[r]&0}$ is a
minimal $\add~M$-resolution of X. If N=DA and T is a convex tilting
of module, then $K\in\add~ T$.}

\vskip0.2in

Now, we can prove Theorem 3 promised in introduction.

\vskip0.2in

{\bf Theorem 4.7.}\ {\it Let A be a hereditary algebra of Dynkin
type over an algebraically closed field $k$, and $T$ be a tilting
$A$-module. Let $\ol T= T\oplus\ol P$ be a tilting $A^{(1)}$-module
as above. Then ${\rm rep.dim \ End}_{A^{(1)}}~\ol T\leq 3$.}

\vskip0.1in

{\bf Proof}\ \  Let $(P_1, P_2, f)$ be the additive generator of
$\Mor~({\rm add}\ A)$. We can assume that $P_2=\Hom(DT, DT_2)$ with
$DT_2\in\add~DT$. Then $(P_1, DT_2, f)$ contains all the
indecomposable projective $T_2(A)$-$\modules$ of the form $(P, DT,
id)$.

Let $M=(P_1,DT_2 , f)\oplus(T, DB,
id)\oplus(DA,0,0)\oplus(0,DB\oplus DT,0)\oplus(T,0,0)\oplus(0,B,0)$
with $B=\End_A~T$. Then $M$ is a generator-cogenerator of
$\Mod~\End_{A^{(1)}}~\ol T$.

Let $(X, Y, f)$ be an indecomposable $\End_{A^{(1)}}~\ol T$-module.
It is easy to see that there are three kinds of indecomposable
$\End_{A^{(1)}}~\ol T$-modules.

\vskip0.1in

(1) The first case is $(X, Y, f)$ such that $X\neq0$ and
$Y\neq0\in\mathscr{T}(\dt)$.

We claim that there is no morphism from $(DA,0,0)$ to $(X,Y,f)$.

In fact, we assume by contrary that $(g,0)\neq 0$ is a morphism from
$(DA,0,0)$ to $(X,Y,f)$, then $\mathrm{Im}~g$ is an injective direct
summand of $X$ and also belongs to $\Ker f$. It follows that
$(\mathrm{Im}~g, 0,0)$ is a direct summand of $(X,Y,f)$ which
contradicts with the assumption of  $(X,Y,f)$.

Let $M_1\rightarrow(X,Y,f)\rightarrow0$ be a minimal
$\add~M$-approximation of $(X,Y,f)$. Then
$M_1=(P',DT',f')\oplus(T_1,DB_1,id)\oplus(0,DB_2\oplus
DT'',0)\oplus(T_3,0,0)$. We get an exact sequence
$$\xymatrix{0\ar[r]&\Ker~\pi_1\ar[r]\ar[d]_h&P'\oplus T_1\oplus T_3\ar[r]^-{\pi_1}\ar[d]&X\ar[r]\ar[d]_f&0\\
          0\ar[r]&\Ker~\pi_2\ar[r]&DT'\oplus DB_1\oplus DB_2\oplus DT''\ar[r]^-{\pi_2}&Y\ar[r]&0}$$
We should mention that the lower row in above commutative diagram
should be the sequence under functor $\Hom_B(DT,-)$. But it doesn't
matter since $\Hom_B(DT,-)$ is an equivalence between $\mathscr{T}(DT)$ and
$\mathscr{Y}(DT)$.

Let $\pi_1=(p,t)$ with $p: P'\rightarrow X$ and $t:T_1\oplus
T_2\rightarrow X$ is an $\add~ T$-approximation of X. According to
Lemma 4.5, $\Ker~t\in\add~ T$.

We consider the following pull-back diagram
$$\xymatrix{0\ar[r]&\Ker~ t\ar[r]\ar@{=}[d]&\Ker~ \pi_1\ar[d]\ar[r]^q&P'\ar[d]_-p\\
            0\ar[r]&\Ker~t\ar[r]&T_1\oplus T_2\ar[r]^-t&X}$$
Note that $\mathrm{Im}~q=P''$ is projective since $A$ is hereditary.
Let $\Ker~t=T_4$. Then $\Ker~\pi_1$=$P''\oplus T_4$.

We claim that $\xymatrix{DT'\oplus DB_1\oplus DB_2\oplus
DT''\ar[r]^-{\pi_2}&Y}$ is an minimal $\add~( DB\oplus
DT)$-approximation of Y and $\pi_2$ is epimorphism.

The reason is that $M_1\rightarrow(X,Y,f)\rightarrow0$ is minimal
and $(0,DB\oplus DT,0)$ is an direct summand of $M$. It follows that
$\pi_2$ is an minimal $\add~( DB\oplus DT)$-approximation of Y. Note
that $Y\in \mathscr{T}(DT)$ which implies that $Y$ is generated by $DT$,
hence $\pi_2$ is epimorphism. By using Lemma 4.6, we know that
$\Ker~\pi_2\in\add~ DT$.

Let $\Ker~\pi_2=DT_5$. Then $DT_5\subset DT'\oplus DT''$, hence
$(\Ker~\pi_1, \Ker~\pi_2, h)$=$(P'', DT_5, h)\oplus(T_4, 0,0)$
belongs to $\add~M$.

\vskip0.1in

(2) The second case is $(X,Y,f)=(X,0,0)$ with $X\in{\rm ind}\ A$.
If $X\in\add~ DA$, then
$0\rightarrow(X,0,0)\rightarrow(X,0,0)\rightarrow0$ is an
$\add~M$-resolution of $(X,0,0)$. If X is not injective, then by
using the same proof as in (1), we can obtain an $\add~M$-resolution
of $(X,0,0)$ with the length at most one.

\vskip0.1in

(3) The third case is $(0,Y,0)$ with $Y\in{\rm ind}\ B$. According
to \cite{APT}, we know that ${\rm rep.dim}~ B\leq 3$ and $B\oplus
DB\oplus DT$ is an Auslander generator for $\Mod\ B$.

Assume that $0\rightarrow N_2\rightarrow N_1\rightarrow Y\rightarrow
0$ is an $\add~ (B\oplus DB\oplus DT)$-resolution of $Y$. It is easy
to check that there is no non-zero morphism from $(P_1,DT_2 ,
f)\oplus(T, DB, id)\oplus(DA,0,0)$ to $(0,Y,0)$, hence
$0\rightarrow(0,N_2,0)\rightarrow(0,N_1,0)\rightarrow(0,Y,0)\rightarrow
0$ is an $\add~M$-resolution of $(0,Y,0)$.

\vskip0.1in

Summary the above discussions, we have shown that any indecomposable
$\End_{A^{(1)}}~\ol T$-module admits an $\add~M$-resolution with the
length at most one. Hence ${\rm rep.dim\ End}_{A^{(1)}}~\ol T\leq
3$. This completes the proof. $\hfill\Box$


\vskip0.2in

\begin{thebibliography}{100}

\bibitem{ABSG} I.Assem, T.Br\"ustle, R.Schiffer, G.Todorov,  Cluster categories and duplicated algebras.
J.Algebra, 305(2006), 548-561.

\bibitem{APT} I.Assem, M.Platzeck, S.Trepode, On the representation dimension of tilted and Laura algebras.
J.Algebra, 296(2006), 426-439.

\bibitem{ASS}  I.Assem, D.Simson, A.Skowronski, Elements of the reprensentation theory of associative
algebras. In: Techniques of Representation Theory, vol.1, Cambridge
University Press, Cambridge, 2006.

\bibitem{Auslander} M.Auslander, The representation dimension of artin algebras.
Queen Mary College Mathematics Notes (1971). Republished in {\it
Selected works of Maurice Auslander.} Amer.Math.Soc., Providence
1999.

\bibitem{AR} M.Auslander, I.Reiten,  Applications of contravariantly finite
subcategories.  Adv. Math., 86(1991), 111-152.

\bibitem{ARS}  M.Auslander, I.Reiten, S.Smal\o, Representation Theory of Artin
Algebras. Cambridge Stud.Adv. Math., vol.36, Cambridge University
Press, 1995.

\bibitem{BM} A.Buan, R.Marsh, M.Reineke, I.Reiten, G.Todorov,
   Tilting and cluster combinatorics. Adv.Math.,204(2006),572-618.


\bibitem{EHIS} K.Erdmann, T.Holm, O.Iyama, J.Schr\"oer, Radical embeddings and representation dimension.
Adv.Math., 185(2004), 159-177.

\bibitem{FGR} R.Fossum, P.Griffith, I.Reiten, Trivival Extension of Abelian
Categories. Lecture Notes in Mathematics 456, Springer-Verlag,
Berlin/New York, 1975.

\bibitem{HR} D.Happel, C.R.Ringel. Tilted algebras. Trans.Amer.Math.Soc., 274(1982), 399-443.

\bibitem{HU} D.Happel, L.Unger, On the quiver of tiltiing modules. J.Algebra, 284(2005), 857-868.

\bibitem{HS} P.Hilton, U.Stammbach, A coures in homological algebra,
Graduate Texts in Mathematices 4, Springer-Verlag, Berlin/New York,
1970.

\bibitem{IT} K.Igusa,G.Todorav. On the finitistic global dimension conjecture for Artin
algebras. Fields Institute Communications, 45(2005), 201-204.

\bibitem{Iyama} O.Iyama. Finiteness of representation dimension. Proc.Amer.Math.Soc., 131(2003), 1011-1014.

\bibitem{llz} X.Lei, H.Lv, S.Zhang, Complements to the almost complete tilting $A^{m}$-modules. Comm.Algebra, 37(2009), 1719-1728.

\bibitem{zl} H.Lv, S.Zhang, Representation dimension of $m$-replicated algebras. Sci.China Math., 53(2010), 1603-1608.

\bibitem{Oppermann} S.Oppermann, Representation dimension of qusi-tilted algebras. J.London Math. Soc., (2)81(2010), 435-456.

\bibitem{Rouquier} R.Rouquier, Representation dimension of exterior algebras. Inv.Math., 165(2006), 357-367.

\bibitem{Xi} C.Xi, On the representation dimension of finite dimensional algebras. J.Algebra, 226(2000), 332-346.

\bibitem{z1} S.Zhang, Tilting mutation and duplicated algebras. Comm.Algebra, 37(2009), 3516-3524.

\bibitem{z2} S.Zhang, Partial tilting modules over m-replicated algebras. J.Algebra, 323(2010), 2538-2546.




\end{thebibliography}
\end{document}